\documentclass[12pt]{article}
\usepackage{amssymb}
\usepackage{mathrsfs}
\usepackage{pifont}
\usepackage{tipa}
\usepackage{amsmath}
\usepackage{amssymb,amsmath}
\usepackage{amsfonts}
\catcode`@=11
 \long\def\@makefntext#1{\noindent #1}

\begin{document}
\title{\large{\textbf{Whittaker modules
for the Schr\"{o}dinger-Virasoro\\ algebra}} }
\author{Xiufu Zhang $^{1,2}$, Shaobin Tan $^1$, Haifeng Lian $^3$\\
{\scriptsize 1. School of Mathematical Sciences, Xiamen University,
Xiamen 361005,  China}\\
{\scriptsize 2. School of Mathematical Sciences, Xuzhou Normal
University, Xuzhou 221116, China}\\
{\scriptsize 3. School of  Computer and Information, Fujian
Agriculture and Forestry University, Fuzhou 350002, China}}
\date{}
\maketitle \footnotetext{\footnotesize* Supported by the National
Natural Science Foundation of China (No. 10671160).}
\footnotetext{\footnotesize** Email: xfzhang@xznu.edu.cn;\quad
tans@xmu.edu.cn;\quad hlian@fjau.edu.cn}

\numberwithin{equation}{section}

\begin{abstract}
In this paper, Whittaker modules for the Schr\"{o}dinger-Virasoro
algebra $\mathfrak{sv}$ are defined. The Whittaker vectors and the
irreducibility of the Whittaker modules are studied. $\mathfrak{sv}$
has a triangular decomposition according to the Cartan algebra
$\mathfrak{h}:$
$$\mathfrak{sv}=\mathfrak{sv}^{-}\oplus\mathfrak{h}\oplus\mathfrak{sv}^{+}.$$
For any Lie algebra homomorphism
$\psi:\mathfrak{sv}^{+}\rightarrow\mathbb{C}$, we can define
Whittaker modules of type $\psi.$ When $\psi$ is nonsingular, the
Whittaker vectors , the irreducibility and the classification of
Whittaker modules are completely determined. When $\psi$ is
singular, by constructing some special Whittaker vectors, we find
that the Whittaker modules are all reducible. Moreover, we get some
more precise results for special $\psi$. \vspace{2mm}\\{\bf 2000
Mathematics Subject Classification:} 17B10, 17B35, 17B65, 17B68
\vspace{2mm}
\\ {\bf Keywords:}  Schr\"{o}dinger-Virasoro algebra,
 Whittaker vector, Whittaker
module, induced module, irreducible module.
\end{abstract}

\vskip 3mm \noindent{\section{{Introduction}}}

\vskip 3mm The Schr\"{o}dinger-Virasoro algebra $\mathfrak{sv},$
playing important roles in mathematics and statistical physics, is a
infinite-dimensional Lie algebra first introduced by M. Henkel in
[7] by looking at the invariance of the free Schr\"{o}dinger
equation. This infinite-dimensional Lie algebra contains both the
Lie algebra of invariance of the free Schr\"{o}dinger equation and
the centerless Virasoro algebra (Witt algebra) as subalgebras. As
natural deformations of the Schr\"{o}dinger-Virasoro algebra
$\mathfrak{sv},$ the twisted Schr\"{o}dinger-Virasoro algebra,
$\varepsilon$-deformation Schr\"{o}dinger-Virasoro algebra,  the
extended Schr\"{o}dinger-Virasoro algebra and  the generalized
Schr\"{o}dinger-Virasoro algebras  are introduced in [20]-[22]. The
derivations, the 2-cocycles, the central extensions and the
automorphisms for these algebras have been well studied by many
authors (e.g., [6],[11], [20]-[23]).

With respect to the representation theory for
Schr\"{o}dinger-Virasoro algebra, the weight modules is well studied
in [10], there it is proved that an irreducible weight module with
finite-dimensional weight spaces over the Schr\"{o}dinger-Virasoro
algebras is a highest/lowest weight module or a uniformly bounded
module. This is the analogue of a well known classical result in the
Virasoro algebra setting conjectured by V. Kac and proved or
partially proved by many authors (see [12], [13] and [19]).

In this paper, we construct and study the so called Whittaker
modules for the Schr\"{o}dinger-Virasoro algebra $\mathfrak{sv}$
which are not weight modules.

The notion of Whittaker modules is first introduced by D. Arnal and
G. Pinczon in [1] in the process of construction of a very vast
family of representations for $sl(2).$ The versions of Whittaker
modules of the complex semisimple Lie algebras are generalized by
Kostant in [9].  The prominent role played by Whittaker modules is
illustrated by the main result in [3] about the classification of
the irreducible modules for $sl_2(\mathbb{C}).$ The result
illustrate that the irreducible $sl_2(\mathbb{C})$-modules fall into
three families: highest (lowest) weight modules, Whittaker modules,
and a third family obtained by localization. Since the construction
of Whittaker modules depends on the triangular decomposition of a
finite-dimensional complex semisimple Lie algebras, it is natural to
consider Whittaker modules for other algebras with a triangular
decomposition. Recently, the Whittaker modules for Virasoro
algebras, Heisenberg algebras, affine Lie algebras as well as
generalized Weyl algebras are studied by M. Ondrus, E. Wiesner, K.
Christodoulopoulou, G. Benkart,  etc.(see [2], [5], [14], [16] and
[17]).

 The Schr\"{o}dinger-Virasoro algebra $\mathfrak{sv}$ has a
triangular decomposition:
$\mathfrak{sv}=\mathfrak{sv}^{-}\oplus\mathfrak{h}\oplus\mathfrak{sv}^{+}$.
For any Lie algebra homomorphism $\psi: \mathfrak{sv}^{+}\rightarrow
\mathbb{C},$ we can define Whittaker modules of type $\psi$ for
$\mathfrak{sv}.$ Moreover, for $\xi\in \mathbb{C},$  we can
construct two special Whittaker modules $W_{\psi}$ and
$L_{\psi,\xi}$ for $\mathfrak{sv}$(see section 2). In section 3 and
section 4 we will study the Whittaker modules of nonsingular type.
In section 3, the Whittaker vectors of $W_{\psi}$ and $L_{\psi,\xi}$
are studied. In section 4, The classification of the irreducible
Whittaker modules of nonsingular type is studied. In the final
section, we study the Whittaker modules of singular type. The
Whittaker vectors of $W_{\psi}$ and $L_{\psi,\xi}$ are studied. By
constructing some special Whittaker vectors, we see that
$L_{\psi,\xi}$ are all reducible. We also get some more precise
results for special $\psi.$

\vskip 3mm Throughout this paper the symbols $\mathbb{C},
\mathbb{N}$, $\mathbb{Z},$ $\mathbb{Z}_{+}$ and $\sum$ represent for
the complex field, the set of nonnegative integers, the set of
integers, the set of positive integers and the sum with finite
summands respectively.

\vskip 3mm\noindent{\section{Definitions and Notations}}

\vskip 3mm The Schr\"{o}dinger-Virasoro algebra $\mathfrak{sv}$ is
defined to be a Lie algebra with $\mathbb{C}$-basis
$\{L_n,M_n,Y_{n+\frac{1}{2}}\mid n\in \mathbb{Z}\}$ subject to the
following Lie brackets:
\begin{align*}
&[L_m,L_n]=(n-m)L_{n+m},\
[L_m,M_n]=nM_{n+m},\\
&[L_m,Y_{n+\frac{1}{2}}]=(n+\frac{1-m}{2})Y_{m+n+\frac{1}{2}},\\
&[Y_{m+\frac{1}{2}},Y_{n+\frac{1}{2}}]=(n-m)M_{m+n+1},\\
&[M_m,M_n]=[M_m,Y_{n+\frac{1}{2}}]=0.
\end{align*}
 It is easy to see the following facts
about $\mathfrak{sv}:$

\vskip 3mm (i) The center of $\mathfrak{sv}$ is $\mathbb{C}M_0.$

\vskip 3mm (ii) $\mathfrak{sv}$ is a semi-direct product of the Witt
algebra $\mathfrak{Vir_0}=\bigoplus_{n\in\mathbb{Z}}\mathbb{C}L_n$
and the two-step nilpotent infinite-dimensional Lie algebra
$\mathfrak{g}=\bigoplus_{n\in\mathbb{Z}}\mathbb{C}M_n\oplus
\bigoplus_{n\in\mathbb{Z}}\mathbb{C}Y_{\frac{1}{2}+n}.$

\vskip 3mm (iii) $\mathfrak{sv}$ has a triangular decomposition
according to the Cartan algebra $\mathfrak{h}=\mathbb{C}L_0\oplus
\mathbb{C}M_0:$
$$\mathfrak{sv}=\mathfrak{sv}^{-}\oplus\mathfrak{h}\oplus\mathfrak{sv}^{+},$$
where
$$\mathfrak{sv}^{+}=\mathrm{span}_{\mathbb{C}}
\{L_n,M_{n},Y_{\frac{1}{2}+m}|m\in \mathbb{N}, n\in
\mathbb{Z}_{+}\},$$
$$\mathfrak{sv}^{-}=\mathrm{span}_{\mathbb{C}}\{L_{-n},M_{-n},
Y_{-\frac{1}{2}-m}|m\in \mathbb{N}, n\in \mathbb{Z}_{+}\}.$$

\vskip 3mm (iv) $\mathfrak{sv}^{+}$ (resp. $\mathfrak{sv}^{-}$) is
generated by $L_1, L_2, M_1$ and $Y_{\frac{1}{2}}$ (resp.
 $L_{-1}, L_{-2}, M_{-1}$ and
$Y_{-\frac{1}{2}}$).

\vskip 3mm In the following of this section we give some notations
which will be frequently used to describe the basis of the universal
enveloping algebra $U(\mathfrak{sv})$ and the basis of Whittaker
modules for the Schr\"{o}dinger-Virasoro algebra. Set
$$\mathfrak{b}^{+}=\mathfrak{sv}^{+}\oplus\mathfrak{h},\ \
\mathfrak{b}^{-}=\mathfrak{sv}^{-}\oplus\mathfrak{h}.$$
 Let $\mathbb{C}[M_0]$ be the polynomial
algebra generated by $M_0.$ Obviously, $\mathbb{C}[M_0]$ is
contained in $Z(\mathfrak{sv}),$ the center of $U(\mathfrak{sv}).$

  As in [18], for a non-decreasing sequence of
  positive integers: $0<\mu_1\leq\mu_2\leq\cdots\leq\mu_s$, we call
   $\mu=(\mu_1,\mu_2,\cdots,\mu_s)$ a $partition,$ and
for a non-decreasing sequence of
 non-negative integers: $0\leq\lambda_1\leq\lambda_2\leq\cdots\leq\lambda_r$, we
 call $\widetilde{\lambda}=(\lambda_1,\lambda_2,\cdots,\lambda_r)$
 a $pseudopartition.$
Let $\mathcal{P}$ denote the set of partitions, and let
$\widetilde{\mathcal{P}}$ represent the set of pseudopartitions.
Then $\mathcal{P}\subset\widetilde{\mathcal{P}}.$ For
$\widetilde{\lambda}\in \widetilde{\mathcal{P}},$ we also write
$\widetilde{\lambda}=(0^{\lambda(0)},1^{\lambda(1)},2^{\lambda(2)}\cdots),$
where $\lambda(k)$ is the number of times of $k$ appears in the
pseudopartition and $\lambda(k)=0$ for $k$ sufficiently large. Then
a pseudopartition $\widetilde{\lambda}$ is a partition whenever
$\lambda(0)=0.$ For $\mu=(\mu_1,\mu_2,\cdots,\mu_s)\in\mathcal{P},$
$\widetilde{\lambda}=(\lambda_1,\lambda_2,\cdots\lambda_r)$ and
$\widetilde{\nu}=(\nu_1,\nu_2,\cdots\nu_t)\in
\widetilde{\mathcal{P}},$ we define
$$|\widetilde{\lambda}|=\lambda_1+\lambda_2+\cdots+\lambda_r,$$
$$\frac{1}{2}+\widetilde{\nu}=(\frac{1}{2}+\nu_1,\frac{1}{2}+\nu_2,\cdots,\frac{1}{2}
+\nu_t),$$
$$|\frac{1}{2}+\widetilde{\nu}|=(\frac{1}{2}+\nu_1)
+(\frac{1}{2}+\nu_2)+\cdots+(\frac{1}{2}+\nu_t),$$
$$\#(\widetilde{\lambda})=\lambda(0)+\lambda(1)+\cdots,$$
$$\#(\mu,\widetilde{\nu},\widetilde{\lambda})
=\#(\mu)+\#(\widetilde{\nu})+\#(\widetilde{\lambda}),$$
$$L_{-\widetilde{\lambda}}=L_{-\lambda_r}\cdots L_{-\lambda_2}L_{-\lambda_1}
=\cdots L_{-2}^{\lambda(2)}L_{-1}^{\lambda(1)}L_{0}^{\lambda(0)},$$
$$M_{-\mu}=M_{-\mu_s}\cdots M_{-\mu_2}M_{-\mu_1}=\cdots M_{-2}^{\mu(2)}M_{-1}^{\mu(1)},$$
$$Y_{-\frac{1}{2}-\widetilde{\nu}}=Y_{-\frac{1}{2}-\nu_t}\cdots Y_{-\frac{1}{2}-\nu_2}
Y_{-\frac{1}{2}-\nu_1}= \cdots
Y_{-\frac{1}{2}-1}^{\nu(1)}Y_{-\frac{1}{2}}^{\nu(0)}.$$ For the sake
of convenience, we define $\bar{0}=(0^{0},1^{0},2^{0},\cdots) $ and
set $L_{\bar{0}}=M_{\bar{0}}=Y_{\frac{1}{2}+\bar{0}}=1\in
U(\mathfrak{sv}).$ In the following, we regard $\bar{0}$ as an
element of $\mathcal{P}$ and $\widetilde{\mathcal{P}}$.

\vskip 3mm For any $(\mu,\widetilde{\nu},\widetilde{\lambda})\in
\mathcal{P}\times\widetilde{\mathcal{P}}\times\widetilde{\mathcal{P}}$
and
$p_{\mu,\widetilde{\nu},\widetilde{\lambda}}(M_0)\in\mathbb{C}[M_0],$
it is obvious that
$$p_{\mu,\widetilde{\nu},\widetilde{\lambda}}(M_0)
M_{-\mu}Y_{-\frac{1}{2}-\widetilde{\nu}}L_{-\widetilde{\lambda}}\in
U(\mathfrak{sv})_{-(|\mu|+|\frac{1}{2}+\widetilde{\nu}|+|\widetilde{\lambda}|)},$$
where $U(\mathfrak{sv})_{a}=\{x\in U(\mathfrak{sv})|[L_0,x]=ax\}$ is
the $a$-weight space of $U(\mathfrak{sv}).$

\vskip 5mm\noindent{\bf{Definition 2.1. }} Let $V$ be a
$\mathfrak{sv}$-module and let
$\psi:\mathfrak{sv}^{+}\rightarrow\mathbb{C}$ be a Lie algebra
homomorphism. A vector $v\in V$ is called a Whittaker vector if
$xv=\psi(x)v$ for every $x\in \mathfrak{sv}^{+}.$  A
$\mathfrak{sv}$-module $V$ is called a Whittaker module of type
$\psi$ if there is a Whittaker vector $w\in V$ which generates $V.$
In this case we call $w$ the cyclic Whittaker vector.

The Lie algebra homomorphism $\psi$ is called nonsingular if
$\psi(M_1)$ is nonzero, otherwise $\psi$ is called singular. The Lie
brackets in the definition of $\mathfrak{sv}$ force
$\psi(L_n)=\psi(M_m)=\psi(Y_{\frac{1}{2}+k})=0$ for $n\geq3, m\geq2,
k\geq1.$

\vskip 3mm For a Lie algebra homomorphism $\psi:
\mathfrak{sv}^{+}\rightarrow\mathbb{C},$ we define
$\mathbb{C}_{\psi}$ to be the one-dimensional
$\mathfrak{sv}^{+}$-module given by $x\alpha=\psi(x)\alpha$ for
$x\in \mathfrak{sv}^{+}$ and $\alpha\in \mathbb{C}.$ Then we have an
induced $\mathfrak{sv}$-module
\begin{equation}
W_{\psi}=U(\mathfrak{sv})\otimes_{U(\mathfrak{sv}^{+})}\mathbb{C}_{\psi}.
\end{equation}
For $\xi\in \mathbb{C},$ $(M_0-\xi)W_{\psi}$ is a submodule of
$W_{\psi}$ since $M_0$ is in the center of $\mathfrak{sv}.$ Set
\begin{eqnarray}
L_{\psi,\xi}:=W_{\psi}/(M_0-\xi)W_{\psi}.
\end{eqnarray}
Then
$L_{\psi,\xi}$ is a quotient module for $\mathfrak{sv}.$ The
following facts about $W_{\psi}$ are obvious:

\vskip 3mm (i) $W_{\psi}$ is a Whittaker module of type $\psi,$ with
cyclic Whittaker vector $w:=1\otimes1;$

\vskip 3mm  (ii) The set
\begin{eqnarray}
\{M_0^{k}M_{-\mu}Y_{-\frac{1}{2}-\widetilde{\nu}}L_{-\widetilde{\lambda}}w|
(\mu, \widetilde{\nu}, \widetilde{\lambda})\in \mathcal{P}\times
\widetilde{\mathcal{P}}\times \widetilde{\mathcal{P}}, k\in
\mathbb{N}\}
\end{eqnarray}
forms a basis of $W_{\psi}$. This follows from the PBW theorem and
the fact that
\begin{eqnarray}
\{M_0^{k}M_{-\mu}Y_{-\frac{1}{2}-\widetilde{\nu}}L_{-\widetilde{\lambda}}|
( \mu,\widetilde{\nu},\widetilde{\lambda})\in
\mathcal{P}\times\widetilde{\mathcal{P}}\times\widetilde{\mathcal{P}},
k\in \mathbb{N}\}
\end{eqnarray}
is a basis of $U(\mathfrak{b}^{-});$

\vskip 3mm (iii) $W_{\psi}$ has the universal property in the sense
that for any Whittaker module $V$ of type $\psi$ generated by
$w^{'},$ there is a surjective homomorphism $\varphi:
W_{\psi}\rightarrow V$ such that $uw\mapsto uw^{'}, \forall u\in
U(\mathfrak{b}^{-}).$ Hence we call $W_{\psi}$ the universal
Whittaker module of type $\psi$.

\vskip 3mm For any $0\neq v=\sum
p_{\mu,\widetilde{\nu},\widetilde{\lambda}}(M_0)
M_{-\mu}Y_{-\frac{1}{2}-\widetilde{\nu}}L_{-\widetilde{\lambda}}w\in
W_{\psi},$ we define
$$maxdeg(v):=\max\{|\mu|+|\frac{1}{2}+\widetilde{\nu}|+|\widetilde{\lambda}||
p_{\mu,\widetilde{\nu},\widetilde{\lambda}}(M_0)\neq0\},$$
$$max_{L_0}(v):=\max\{\lambda(0)|
p_{\mu,\widetilde{\nu},\widetilde{\lambda}}(M_0)\neq0\}.$$ We set
$maxdeg(w)=0,\ maxdeg(0)=-\infty.$

\vskip 3mm\noindent{\bf{Remark 2.2.}} For any $x\in
U(\mathfrak{sv}^{+}),$ $w^{'}=uw, u\in U(\mathfrak{b}^{-}),$ we have
$$(x-\psi(x))w^{'}=[x,u]w.$$  In particular,
$$(E_{n}-\psi(E_{n}))w^{'}=[E_n,u]w,$$ where $E_{n}=L_{n}$ or $M_{n}$
or $Y_{\frac{1}{2}+(n-1)}, \forall n\in \mathbb{Z}_{+}.$

\vskip 3mm For $m\in\mathbb{Z}_{+}, n, k\in \mathbb{N},$ $\mu\in
\mathcal{P},$ $\widetilde{\lambda}, \widetilde{\nu}\in
\widetilde{\mathcal{P}},$ we give some identities of
$U(\mathfrak{sv}),$ each of them can be checked by induction on
$\#(\widetilde{\lambda})$ or $\#(\widetilde{\nu})$ or $a\in
\mathbb{N}:$
\begin{equation}
M_{m}L_{-\widetilde{\lambda}}=\sum
a_{i}M_{-m_{i}}L_{-\widetilde{\lambda^{'}}^{(i)}}+\sum
b_{i}L_{-\widetilde{\lambda^{''}}^{(i)}}M_{n_{i}}+L_{-\widetilde{\lambda}}M_{m},
\end{equation}
where $a_{i}, b_{i}\in \mathbb{C},$ $m_{i}\geq0,$ $0<n_{i}\leq m,$
$|\widetilde{\lambda^{'}}^{(i)}|+m_{i}=|\widetilde{\lambda^{''}}^{(i)}|-n_{i}
=|\widetilde{\lambda}|-m,$ and ${\lambda^{''}}^{(i)}(0)<\lambda(0)$
if $n_{i}=m.$
\begin{equation} \label{eq:1}
M_{m}L_{-k}^{a}=\sum_{i=0}^{a}(-1)^{i}(\prod_{j=0}^{i-1}(m-jk))
(^{a}_{i})L_{-k}^{a-i}M_{m-ik}.
\end{equation}
\begin{equation} \label{eq:1}
Y_{\frac{1}{2}+n}L_{-\widetilde{\lambda}}=\sum
a_{i}Y_{-\frac{1}{2}-m_i}L_{-\widetilde{\lambda^{'}}^{(i)}}+\sum
b_{i}L_{-\widetilde{\lambda^{''}}^{(i)}}Y_{\frac{1}{2}+n_{i}}
+L_{-\widetilde{\lambda}}Y_{\frac{1}{2}+n},
\end{equation}
where $a_{i}, b_{i}\in \mathbb{C},$ $0\leq n_{i}\leq n,$
$|\widetilde{\lambda^{'}}^{(i)}|+(\frac{1}{2}+m_{i})
=|\widetilde{\lambda^{''}}^{(i)}|-(\frac{1}{2}+n_{i})=|\widetilde{\lambda}|-(\frac{1}{2}+n),$
and ${\lambda^{''}}^{(i)}(0)<\lambda(0)$ if $n_{i}=n.$
\begin{equation} \label{eq:1}
Y_{\frac{1}{2}+n}Y_{-\frac{1}{2}-\widetilde{\nu}}=\sum
b_{i}Y_{-\frac{1}{2}-\widetilde{\nu^{'}}^{(i)}}M_{n_{i}}+
Y_{-\frac{1}{2}-\widetilde{\nu}}Y_{\frac{1}{2}+n},
\end{equation}
where $b_{i}\in \mathbb{C},$ $n_{i}\leq n,$
$|\frac{1}{2}+\widetilde{\nu^{'}}^{(i)}|-n_{i}=
|\frac{1}{2}+\widetilde{\nu}|-(\frac{1}{2}+n).$

\begin{equation} \label{eq:1}
L_{n}L_{-\widetilde{\lambda}}=\sum
a_{i}L_{-\widetilde{\lambda^{''}}^{(i)}}L_{n_{i}}+L_{-\widetilde{\lambda}}L_{n},
\end{equation}
where $a_{i}\in \mathbb{C},$ $n_{i}\leq n,$
$\widetilde{\lambda^{''}}^{(i)}-n_{i}=|\widetilde{\lambda}|-n,$ and
${\lambda^{''}}^{(i)}(0)<\lambda(0)$ if $n_{i}=n.$

\begin{equation} \label{eq:1}
L_{n}M_{-\mu}=\sum a_{i}M_{-{\mu^{'}}^{(i)}}+\sum
b_{i}M_{-{\mu^{''}}^{(i)}}M_{m_{i}}+M_{-\mu}L_{n},
\end{equation}
where $a_{i}, b_{i}\in \mathbb{C},$ $m_{i}<n,$
$|{\mu^{'}}^{(i)}|=|{\mu^{''}}^{(i)}|-m_{i}=|\mu|-n.$

\begin{equation} \label{eq:1}
L_{n}Y_{-\frac{1}{2}-\widetilde{\nu}}=\sum
a_{i}Y_{-\frac{1}{2}-\widetilde{\nu^{'}}^{(i)}}+\sum
b_{i}Y_{-\frac{1}{2}-{\widetilde{\nu^{''}}}^{(i)}}Y_{\frac{1}{2}+n_{i}}+\sum
c_{i}Y_{-\frac{1}{2}-{\widetilde{\nu^{''}}}^{(i)}}M_{m_i}
+Y_{-\frac{1}{2}-\widetilde{\nu}}L_{n},
\end{equation}
where $a_{i}, b_{i}, c_{i}\in \mathbb{C},$ $m_i, n_{i}<n,
|\frac{1}{2}+\widetilde{\nu^{'}}^{(i)}|
=|\frac{1}{2}+\widetilde{\nu^{''}}^{(i)}|-(\frac{1}{2}+n_i)
=|\frac{1}{2}+\widetilde{\nu}|-n.$

\vskip 5mm\noindent{\section{Whittaker vectors for Whittaker modules
of nonsingular type}}

In this section we always assume that the Lie homomorphism $\psi$ is
nonsingular, that is $\psi(M_1)\neq 0.$ Let $W_{\psi}$ and
$L_{\psi,\xi}$ be the Whittaker modules for Schr\"{o}dinger-Virasoro
$\mathfrak{sv}$ defined by (2.1) and (2.2) respectively. The main
results of this section are given in Theorem 3.5 and Theorem 3.7 in
which we characterize the Whittaker vectors in $W_{\psi}$ and
$L_{\psi,\xi}.$ For this purpose, we first give a series lemmas
which will be used to prove our main results.

\vskip 3mm\noindent{\bf{Lemma 3.1.}} Let $E_n$ be defined in Remark
2.2, $w=1\otimes1\in W_{\psi}$ be the cyclic Whittaker vector. For
$n\in \mathbb{Z}_{+},$
$$E_nM_{-\mu}Y_{-\frac{1}{2}-\widetilde{\nu}}L_{-\widetilde{\lambda}}w
=v^{'}+v^{''}+\psi(E_n)M_{-\mu}
Y_{-\frac{1}{2}-\widetilde{\nu}}L_{-\widetilde{\lambda}}w,$$ where
$maxdeg(v^{'})<|\mu|+|\frac{1}{2}+\widetilde{\nu}|+|\widetilde{\lambda}|,$
$max_{L_0}(v^{''})<\lambda(0).$

\vskip 3mm\noindent{\bf{Proof.}} If $E_{n}=M_{n},$ the result
follows from (2.5). If $E_{n}=Y_{\frac{1}{2}+(n-1)},$ it follows
from (2.8), (2.5) and (2.7). If $E_{n}=L_{n},$ it follows from
(2.10), (2.11), (2.9), (2.5) and (2.7).\hfill$\Box$

\vskip 3mm\noindent{\bf{Lemma 3.2.}} (i) For $m\in \mathbb{Z}_{+},
\widetilde{\lambda}\in \widetilde{\mathcal{P}},$ then
$maxdeg(M_mL_{-\widetilde{\lambda}}w)
\leq|\widetilde{\lambda}|-m+1;$

\vskip 3mm (ii) For $a, k\in \mathbb{N},$ then
$$[M_{k+1},L_{-k}^{a}]w=v-a(k+1)\psi(M_1)L_{-k}^{a-1}w,$$ where
$maxdeg(v)<(a-1)k$ if $k>0,$ and $max_{L_0}(v)<a-1$ if $k=0;$

\vskip 3mm (iii) Suppose
$\widetilde{\lambda}=(k^{\lambda(k)},(k+1)^{\lambda(k+1)},\cdots),$
$\lambda(k)\neq0.$ Then $$[M_{k+1},L_{-\widetilde{\lambda}}]w=
v-(k+1)\lambda(k)\psi(M_1)L_{-\widetilde{\lambda^{'}}}w,$$ where
$\widetilde{\lambda^{'}}$ satisfies $\lambda^{'}(k)=\lambda(k)-1,
\lambda^{'}(i)=\lambda(i)$ for all $i>k,$
$maxdeg(v)<|\widetilde{\lambda}|-k$ if $k>0$ or $v=v^{'}+v^{''}$
with $maxdeg(v^{'})<|\widetilde{\lambda}|-k$ and
 $max_{L_0}(v^{''})<\lambda(k)-1$ if $k=0;$

\vskip 3mm\noindent{\bf{Proof.}} (i) follows from (2.5) and the fact
that $\psi(M_{i})=0$ if $i\geq2.$ (ii) follows from (2.6). For
(iii), we denote
$L_{-\widetilde{\lambda}}=L_{-\widetilde{\lambda^{'}}}L_{-k}^{\lambda(k)}.$
Then
\begin{equation} \label{eq:1}
[M_{k+1},L_{-\widetilde{\lambda}}]w
=[M_{k+1},L_{-\widetilde{\lambda^{'}}}]L_{-k}^{\lambda(k)}w+
L_{-\widetilde{\lambda^{'}}}[M_{k+1},L_{-k}^{\lambda(k)}]w.
\end{equation}
By using the assumption of $k$, we see that
$[M_{k+1},L_{-\widetilde{\lambda^{'}}}]\in U(\mathfrak{b}^{-})$ and
$$maxdeg([M_{k+1},L_{-\widetilde{\lambda^{'}}}]L_{-k}^{\lambda(k)}w)<|\widetilde{\lambda}|-k.$$
For the second term on the right hand side of (3.1), by using (ii)
we see that
$$L_{-\widetilde{\lambda^{'}}}[M_{k+1},L_{-k}^{\lambda(k)}]w
=L_{-\widetilde{\lambda^{'}}}v-a(k+1)\psi(M_1)L_{-\widetilde{\lambda^{'}}}L_{-k}^{a-1}w,$$
where
$$maxdeg(L_{-\widetilde{\lambda^{'}}}v)<(a-1)k
+|\widetilde{\lambda^{'}}|=|\widetilde{\lambda}|-k$$ if $k>0,$ and
$$max_{L_0}(L_{-\widetilde{\lambda^{'}}}v)<a-1=\lambda(0)-1$$
if $k=0.$ Thus (iii) holds.\hfill$\Box$

\vskip 3mm \noindent{\bf{Lemma 3.3.}} For $m, k\in\mathbb{N},
\tilde{\nu},\tilde{\lambda}\in \tilde{\mathcal{P}},$ we have

 (i) $maxdeg([Y_{\frac{1}{2}+m},
Y_{-\frac{1}{2}-\widetilde{\nu}}L_{-\widetilde{\lambda}}]w)\leq
|\frac{1}{2}+\widetilde{\nu}|+|\widetilde{\lambda}|-(\frac{1}{2}+m)+1;$

\vskip 3mm (ii) If $\nu(i)=\lambda(i)=0$ for all $0\leq i\leq k,$
then
$$maxdeg([Y_{\frac{1}{2}+k+1},
Y_{-\frac{1}{2}-\widetilde{\nu}}L_{-\widetilde{\lambda}}]w)
\leq|\frac{1}{2}+\widetilde{\nu}|+|\widetilde{\lambda}|-k-1;$$

\vskip 3mm (iii) If $\lambda(i)=0$ for all $0\leq i\leq k,$
$\nu(j)=0$ for all $0\leq j<k$ and $\nu(k)\neq0,$ then
$$[Y_{\frac{1}{2}+k+1},
Y_{-\frac{1}{2}-\widetilde{\nu}}L_{-\widetilde{\lambda}}]w
=v-2(1+k)\psi(M_{1})\nu(k)
Y_{-\frac{1}{2}-\widetilde{\nu^{'}}}L_{-\widetilde{\lambda}},$$
where $maxdeg(v)<|\widetilde{\lambda}|+
|\frac{1}{2}+\widetilde{\nu}|-\frac{1}{2}-k$, $\nu^{'}$ satisfies
that $\nu^{'}(i)=\nu(i)$ for all $i\neq k$ and
$\nu^{'}(k)=\nu(k)-1.$

\vskip 3mm\vskip 3mm\noindent{\bf{Proof.}} For (i), note that
\begin{equation} \label{eq:1}
[Y_{\frac{1}{2}+m},
Y_{-\frac{1}{2}-\widetilde{\nu}}L_{-\widetilde{\lambda}}]w=[Y_{\frac{1}{2}+m},
Y_{-\frac{1}{2}-\widetilde{\nu}}]L_{-\widetilde{\lambda}}w+
Y_{-\frac{1}{2}-\widetilde{\nu}}[Y_{\frac{1}{2}+m},L_{-\widetilde{\lambda}}]w.
\end{equation}
By using (2.8) and Lemma 3.2 (i) to the first term on the right hand
side of (3.2), we see that
$$maxdeg([Y_{\frac{1}{2}+m},
Y_{-\frac{1}{2}-\widetilde{\nu}}]L_{-\widetilde{\lambda}}w) \leq
|\frac{1}{2}+\widetilde{\nu}|+|\widetilde{\lambda}|-(\frac{1}{2}+m)+1.$$
By using (2.7) to the second term on the right hand side of (3.2),
we see that
$$maxdeg(
Y_{-\frac{1}{2}-\widetilde{\nu}}[Y_{\frac{1}{2}+m},L_{-\widetilde{\lambda}}]w)
\leq
|\frac{1}{2}+\widetilde{\nu}|+|\widetilde{\lambda}|-(\frac{1}{2}+m)+\frac{1}{2}.$$
Thus (i) holds.

\vskip 3mm For (ii), by using the assumption of $k$, we see that
$[Y_{\frac{1}{2}+k+1},Y_{-\frac{1}{2}-\widetilde{\nu}}]\in
U(\mathfrak{b}^{-}).$ Thus
$$maxdeg([Y_{\frac{1}{2}+k+1},Y_{-\frac{1}{2}-\widetilde{\nu}}]
L_{-\widetilde{\lambda}}w)\leq|\widetilde{\lambda}|+|\frac{1}{2}+\widetilde{\nu}|
-(\frac{1}{2}+k+1).$$ By using (2.7), we see that
$$maxdeg(Y_{-\frac{1}{2}-\widetilde{\nu}}
[Y_{\frac{1}{2}+k+1},L_{-\widetilde{\lambda}}]w)
\leq|\widetilde{\lambda}|+|\frac{1}{2}+\widetilde{\nu}|-(k+1).$$
Thus (ii) follows.

\vskip 3mm Finally, for (iii), we denote
$Y_{-\frac{1}{2}-\widetilde{\nu}}=Y_{-\frac{1}{2}-\widetilde{\nu^{''}}}
Y_{-\frac{1}{2}-k}^{\nu(k)}.$ Then
\begin{eqnarray*}[Y_{\frac{1}{2}+k+1},
Y_{-\frac{1}{2}-\widetilde{\nu}}L_{-\widetilde{\lambda}}]w
=[Y_{\frac{1}{2}+k+1},
Y_{-\frac{1}{2}-\widetilde{\nu^{''}}}]Y_{-\frac{1}{2}-k}^{\nu(k)}
L_{-\widetilde{\lambda}}w+
Y_{-\frac{1}{2}-\widetilde{\nu^{''}}}[Y_{\frac{1}{2}+k+1},Y_{-\frac{1}{2}-k}^{\nu(k)}]
L_{-\widetilde{\lambda}}w
\end{eqnarray*}
\begin{equation} \label{eq:1}+
Y_{-\frac{1}{2}-\widetilde{\nu^{''}}}Y_{-\frac{1}{2}-k}^{\nu(k)}[Y_{\frac{1}{2}+k+1},
L_{-\widetilde{\lambda}}]w.
\end{equation}
For the first term on the right hand side of (3.3), since
$[Y_{\frac{1}{2}+k+1}, Y_{-\frac{1}{2}-\widetilde{\nu^{''}}}]\in
U(\mathfrak{b}^{-}),$ we see that
$$maxdeg([Y_{\frac{1}{2}+k+1},
Y_{-\frac{1}{2}-\widetilde{\nu^{''}}}]Y_{-\frac{1}{2}-k}^{\nu(k)}
L_{-\widetilde{\lambda}}w)\leq|\frac{1}{2}+\widetilde{\nu}|
+|\widetilde{\lambda}|-(\frac{1}{2}+k+1).$$ For the second term,
since $[Y_{\frac{1}{2}+k+1},Y_{-\frac{1}{2}-k}^{\nu(k)}]
=-2(k+1)Y_{-\frac{1}{2}-k}^{\nu(k)-1}M_{1}$ and
$[M_{1},L_{-\widetilde{\lambda}}]\in U(\mathfrak{b}^{-})$ according
to the assumption of $k,$ by using (2.5), we see that
$$Y_{-\frac{1}{2}-\widetilde{\nu^{''}}}[Y_{\frac{1}{2}+k+1},Y_{-\frac{1}{2}-k}^{\nu(k)}]
L_{-\widetilde{\lambda}}w=v^{'}-2(k+1)\psi(M_{1})\nu(k)
Y_{-\frac{1}{2}-\widetilde{\nu^{'}}}L_{-\widetilde{\lambda}}w,$$
where
$maxdeg(v^{'})\leq|\widetilde{\lambda}|+|\frac{1}{2}+\widetilde{\nu}|
-(\frac{1}{2}+k+1),$ $\nu^{'}$ satisfies that $\nu^{'}(i)=\nu(i)$
for all $i\neq k$ and $\nu^{'}(k)=\nu(k)-1.$

\noindent For the third term, note that
$\psi(Y_{\frac{1}{2}+k+1})=0,$ by using (2.7), we see that
$$maxdeg(Y_{-\frac{1}{2}-\widetilde{\nu^{''}}}Y_{-\frac{1}{2}-k}^{\nu(k)}[Y_{\frac{1}{2}+k+1},
L_{-\widetilde{\lambda}}]w)\leq
|\frac{1}{2}+\widetilde{\nu}|+|\widetilde{\lambda}|-k-1.$$ Thus
(iii) follows.\hfill$\Box$

\vskip 3mm\noindent{\bf{Lemma 3.4.}} For $m\in \mathbb{N},\mu\in
\mathcal{P}, \tilde{\nu}\in \tilde{\mathcal{P}},$ we have
$$maxdeg([L_m,M_{-\mu} Y_{-\frac{1}{2}-\widetilde{\nu}}]w)
\leq|\mu|+|\frac{1}{2}+\widetilde{\nu}|-m+1.$$

\vskip 3mm\noindent{\bf{Proof.}}   By (2.10) and (2.11),  we can
write $L_mM_{-\mu} Y_{-\frac{1}{2}-\widetilde{\nu}}$ as a linear
combination of the PBW basis (2.4) of $U(\mathfrak{sv}):$
\begin{eqnarray*}
L_mM_{-\mu}Y_{-\frac{1}{2}-\widetilde{\nu}}&=&
\sum_{\mu^{'},\widetilde{\nu^{'}}}
p_{\mu^{'},\widetilde{\nu^{'}}}(M_0)
M_{-\mu^{'}}Y_{-\frac{1}{2}-\widetilde{\nu^{'}}}
\\
&+&\sum_{\mu^{''},\widetilde{\nu^{''}},n,E_n}
p_{\mu^{''},\widetilde{\nu^{''}},
n}(M_0)M_{-\mu^{''}}Y_{-\frac{1}{2}-\widetilde{\nu^{''}}}E_n,
\end{eqnarray*} where $n\in \mathbb{Z}_{+},$
$\mu^{'},\mu^{''}\in \mathcal{P},\widetilde{\nu^{'}},
\widetilde{\nu^{''}}\in \widetilde{\mathcal{P}}$ satisfying
$|\mu^{'}|+|\frac{1}{2}+\widetilde{\nu^{'}}|=
|\mu^{''}|+|\frac{1}{2}+\widetilde{\nu^{''}}|-n
=|\mu|+|\frac{1}{2}+\widetilde{\nu}|-m;$ $E_n=M_n$ or
$Y_{\frac{1}{2}+(n-1)}.$ Noting that $M_iw=Y_{\frac{1}{2}+j}w=0$ for
$i>1, j>0,$ we see that Lemma 3.4 holds. \hfill$\Box$

\vskip 3mm\noindent{\bf{Theorem 3.5.}} Suppose $\psi(M_1)\neq0$ and
$W_{\psi}$ is the universal Whittaker module for $\mathfrak{sv}$
with cyclic Whittaker vector $w=1\otimes1.$ Then $v\in W_{\psi}$ is
a Whittaker vector if and only if $v=uw$ for some $u\in
\mathbb{C}[M_0].$

\vskip 3mm\noindent{\bf{Proof.}} It is obvious that $uw$ is a
Whittaker vector if $u\in \mathbb{C}[M_0]$ as $M_0$ is in the center
of $\mathfrak{sv}.$

Let $w^{'}\in W_{\psi}$ be an arbitrary vector. We can write $w^{'}$
as a linear combination of the basis (2.3) of $W_{\psi}$:
\begin{equation} \label{eq:1}
w^{'}=\sum_{\mu,\widetilde{\nu},\widetilde{\lambda}}
p_{\mu,\widetilde{\nu},\widetilde{\lambda}}(M_0)
M_{-\mu}Y_{-\frac{1}{2}-\widetilde{\nu}}L_{-\widetilde{\lambda}}w,
\end{equation} where
$p_{\mu,\widetilde{\nu},\widetilde{\lambda}}(M_0)\in
\mathbb{C}[M_0].$ Set
$$N:=max\{|\mu|+|\frac{1}{2}+\widetilde{\nu}|+|\widetilde{\lambda}||
p_{\mu,\widetilde{\nu},\widetilde{\lambda}}(M_0)\neq0\},$$
$$\Lambda_N:=\{(\mu,\widetilde{\nu},\widetilde{\lambda})|
 p_{\mu,\widetilde{\nu},\widetilde{\lambda}}(M_0)\neq0,
|\mu|+|\frac{1}{2}+\widetilde{\nu}|+|\widetilde{\lambda}|=N\}.$$ We
first show that the Whittaker vectors in $W_{\psi}$ are all of type
$\psi.$ In fact, let
$\psi^{'}:\mathfrak{sv}^{+}\rightarrow\mathbb{C}$ be a Lie algebra
homomorphism which is different from $\psi.$ Then there exists at
least one element in $\{L_1, L_2, M_1, Y_{\frac{1}{2}}\},$ denoted
by $E,$ such that $\psi(E)\neq \psi^{'}(E).$ Assume $w^{'}$ is a
Whittaker vector of type $\psi^{'},$ then by the definition we have
$$
Ew^{'}=\psi^{'}(E)w^{'}=\sum_{(\mu,\widetilde{\nu},
\widetilde{\lambda})\notin\Lambda_{N}}
p_{\mu,\widetilde{\nu},\widetilde{\lambda}}(M_0)
M_{-\mu}Y_{-\frac{1}{2}-\widetilde{\nu}}L_{-\widetilde{\lambda}}\psi^{'}(E)w$$
\begin{equation}+\sum_{(\mu,\widetilde{\nu},\widetilde{\lambda})\in\Lambda_{N}}
p_{\mu,\widetilde{\nu},\widetilde{\lambda}}(M_0)
M_{-\mu}Y_{-\frac{1}{2}-\widetilde{\nu}}L_{-\widetilde{\lambda}}\psi^{'}(E)w.
\end{equation}
On the other hand, if we denote
$$K:=max\{\lambda(0)|(\mu,\widetilde{\nu},\widetilde{\lambda})\in\Lambda_{N},
p_{\mu,\widetilde{\nu},\widetilde{\lambda}}(M_{0})\neq0\},$$ then
by Remark 2.2 and Lemma 3.1 we have
\begin{eqnarray}Ew^{'}=v^{'}+v^{''}+\sum_{\begin{subarray}\
(\mu,\widetilde{\nu},\widetilde{\lambda})\in\Lambda_{N}  \\
\ \ \lambda(0)=K
\end{subarray}}
p_{\mu,\widetilde{\nu},\widetilde{\lambda}}(M_0)
M_{-\mu}Y_{-\frac{1}{2}-\widetilde{\nu}}L_{-\widetilde{\lambda}}\psi(E)w,
\end{eqnarray}
where $maxdeg(v^{'})<N,$ $max_{L_0}(v^{''})<K.$  By comparing (3.5)
and (3.6) we obtain $\psi^{'}(E)=\psi(E),$ which is a contradiction
to our assumption that $\psi^{'}(E)\neq\psi(E).$

 Next, for $w^{'}$ defined in (3.4), we want to show that if there is
$(\bar{0},\bar{0},\bar{0})\neq(\mu,\widetilde{\nu},\widetilde{\lambda})\in
\mathcal{P}\times\widetilde{\mathcal{P}}\times\widetilde{\mathcal{P}}$
such that $p_{\mu,\widetilde{\nu},\widetilde{\lambda}}(M_0)\neq 0,$
then there is $E_n\in\{L_n, M_n,
Y_{\frac{1}{2}+(n-1)}|n\in\mathbb{Z}_{+}\}$ such that
$(E_n-\psi(E_n))w^{'}\neq0,$  which will prove the necessity.

 Assume that
$p_{\mu,\widetilde{\nu},\widetilde{\lambda}}(M_0)\neq 0$ for some
$(\mu,\widetilde{\nu},\widetilde{\lambda})\neq
(\bar{0},\bar{0},\bar{0}).$ By Remark 2.2,
\begin{eqnarray*}
(E_n-\psi(E_n))w^{'}= \sum_{\mu,\widetilde{\nu},\widetilde{\lambda}}
p_{\mu,\widetilde{\nu},\widetilde{\lambda}}(M_0) [E_n,
M_{-\mu}Y_{-\frac{1}{2}-\widetilde{\nu}}L_{-\widetilde{\lambda}}]w.
\end{eqnarray*}
Set
$$\underline{k}:=min\{n\in \mathbb{N}|
\mu(n)\neq0 \ \mathrm{or} \ \nu(n)\neq0\ \mathrm{or} \
\lambda(n)\neq0 \ \mathrm{for} \  \mathrm{some} \
(\mu,\widetilde{\nu},\widetilde{\lambda})\in \Lambda_N\}.$$

\noindent We divide the argument into three cases.

\vskip 3mm\noindent{\bf{Case I.}} $\underline{k}$ satisfies
$\lambda(\underline{k})\neq0$ for some
$(\mu,\widetilde{\nu},\widetilde{\lambda})\in \Lambda_N.$

We have
\begin{eqnarray*}
&&(M_{\underline{k}+1}-\psi(M_{\underline{k}+1}))w^{'}\\
&=&\sum_{(\mu,\widetilde{\nu},\widetilde{\lambda})\notin\Lambda_N}
p_{\mu,\widetilde{\nu},\widetilde{\lambda}}
(M_0)[M_{\underline{k}+1},
M_{-\mu}Y_{-\frac{1}{2}-\widetilde{\nu}}L_{-\widetilde{\lambda}}]w\\
&+&\sum_{\begin{subarray}\ (\mu,\widetilde{\nu},\widetilde{\lambda})\in\Lambda_N\\
\ \ \lambda(\underline{k})=0
\end{subarray}} p_{\mu,\widetilde{\nu},\widetilde{\lambda}}
(M_0)[M_{\underline{k}+1},
M_{-\mu}Y_{-\frac{1}{2}-\widetilde{\nu}}L_{-\widetilde{\lambda}}]w\\
\end{eqnarray*}
\begin{equation} \label{eq:1}
+\sum_{\begin{subarray}\ (\mu,\widetilde{\nu},\widetilde{\lambda})\in\Lambda_N\\
\ \ \lambda(\underline{k})\neq0
\end{subarray}} p_{\mu,\widetilde{\nu},\widetilde{\lambda}}
(M_0)[M_{\underline{k}+1},
M_{-\mu}Y_{-\frac{1}{2}-\widetilde{\nu}}L_{-\widetilde{\lambda}}]w
\end{equation} For the first term on the right
hand side of (3.7), by using Lemma 3.2 (i), we know that the degree
of it is strictly smaller than $N-\underline{k}.$ For the second
term on the right hand side of (3.7), note that $\lambda(i)=0$ for
$0\leq i\leq \underline{k},$ we have
$$[M_{\underline{k}+1},M_{-\mu}Y_{-\frac{1}{2}-\widetilde{\nu}}L_{-\widetilde{\lambda}}]
=M_{-\mu}Y_{-\frac{1}{2}-\widetilde{\nu}}[M_{\underline{k}+1},L_{-\widetilde{\lambda}}]
\in U(\mathfrak{b}^{-}).$$ Thus the degree of it is also strictly
smaller than $N-\underline{k}.$ Now using Lemma 3.2 (iii) to the
third term on the right hand side of (3.7), we know that it is of
the form
$$v-\sum_{\begin{subarray}
\ \ (\mu,\widetilde{\nu},\widetilde{\lambda})\in\Lambda_N\\
\ \ \lambda(\underline{k})\neq0
\end{subarray}}(\underline{k}+1)\lambda(\underline{k})\psi(M_1)p_{\mu,\widetilde{\nu},\widetilde{\lambda}}
(M_0)M_{-\mu}Y_{-\frac{1}{2}-\widetilde{\nu}}L_{-\widetilde{\lambda^{'}}}w,$$
where if $\underline{k}=0$ then $v=v^{'}+v^{''}$ such that
$maxdeg(v^{'})<N-\underline{k}$ and
$max_{L_0}(v^{''})<\lambda(0)-1,$ if $\underline{k}>0$ then
$maxdeg(v)<N-\underline{k};$ $\widetilde{\lambda^{'}}$ satisfies
$\lambda^{'}(\underline{k})=\lambda(\underline{k})-1,\lambda^{'}(i)=\lambda(i)$
for all $i>\underline{k}.$ Thus the degree of the third term is
equal to $N-\underline{k}.$ This proves
$(M_{\underline{k}+1}-\psi(M_{\underline{k}+1}))w^{'}\neq0.$

\vskip 3mm\noindent{\bf{Case II.}} $\underline{k}$ satisfies
$\nu(\underline{k})\neq0$ for some
$(\mu,\widetilde{\nu},\widetilde{\lambda})\in \Lambda_N$ and
$\lambda(\underline{k})=0$ for any
$(\mu,\widetilde{\nu},\widetilde{\lambda})\in \Lambda_N.$

In this case, we use
$Y_{\frac{1}{2}+\underline{k}+1}-\psi(Y_{\frac{1}{2}+\underline{k}+1})$
to act on both sides of (3.4), then
\begin{eqnarray*}
&&(Y_{\frac{1}{2}+\underline{k}+1}-\psi(Y_{\frac{1}{2}+\underline{k}+1}))w^{'}\\
&=&\sum_{(\mu,\widetilde{\nu},\widetilde{\lambda})\notin\Lambda_N}
p_{\mu,\widetilde{\nu},\widetilde{\lambda}}
(M_0)[Y_{\frac{1}{2}+\underline{k}+1},
M_{-\mu}Y_{-\frac{1}{2}-\widetilde{\nu}}L_{-\widetilde{\lambda}}]w\\
&+&\sum_{\begin{subarray}
\ (\mu,\widetilde{\nu},\widetilde{\lambda})\in\Lambda_N\\
\ \ \nu(\underline{k})=0
\end{subarray}} p_{\mu,\widetilde{\nu},\widetilde{\lambda}}
(M_0)[Y_{\frac{1}{2}+\underline{k}+1},
M_{-\mu}Y_{-\frac{1}{2}-\widetilde{\nu}}L_{-\widetilde{\lambda}}]w\\
\end{eqnarray*}
\begin{equation} \label{eq:1}
\ +\sum_{\begin{subarray}
\ (\mu,\widetilde{\nu},\widetilde{\lambda})\in\Lambda_N\\
\ \ \nu(\underline{k})\neq0
\end{subarray}} p_{\mu,\widetilde{\nu},\widetilde{\lambda}}
(M_0)[Y_{\frac{1}{2}+\underline{k}+1},
M_{-\mu}Y_{-\frac{1}{2}-\widetilde{\nu}}L_{-\widetilde{\lambda}}]w.
\end{equation}
By using Lemma 3.3 (i) to the first term on the right hand side of
(3.8), Lemma 3.3 (ii) to the second term and Lemma 3.3 (iii) to the
third term, we have
$$(Y_{\frac{1}{2}+\underline{k}+1}-\psi(Y_{\frac{1}{2}+\underline{k}+1}))w^{'}
=v-\sum_{\begin{subarray}
\ (\mu,\widetilde{\nu},\widetilde{\lambda})\in\Lambda_N\\
\ \ \nu(\underline{k})\neq0
\end{subarray}}2(\underline{k}+1)\nu(\underline{k})\psi(M_1)
p_{\mu,\widetilde{\nu},\widetilde{\lambda}}
(M_0)M_{-\mu}Y_{-\frac{1}{2}-\widetilde{\nu^{'}}}L_{-\widetilde{\lambda}}w,$$
where $maxdeg(v)<N-\frac{1}{2}-\underline{k};$ $\nu^{'}(i)=\nu(i)$
for all $i\neq \underline{k}$ and
$\nu^{'}(\underline{k})=\nu(\underline{k})-1.$ Thus
$(Y_{\frac{1}{2}+\underline{k}+1}-\psi(Y_{\frac{1}{2}+\underline{k}+1}))w^{'}\neq0.$

\vskip 3mm\noindent{\bf{Case III.}} $\underline{k}$ satisfies
$\mu(\underline{k})\neq0$ for some
$(\mu,\widetilde{\nu},\widetilde{\lambda})\in \Lambda_N$ and
$\lambda(\underline{k})=\nu(\underline{k})=0$ for any
$(\mu,\widetilde{\nu},\widetilde{\lambda})\in \Lambda_N.$ Note that
in this case $\underline{k}>0$ since $\mu\in \mathcal{P}.$

\vskip 3mm\noindent{\bf{Subcase 1.}} $\widetilde{\lambda}=\bar{0}$
for any $(\mu,\widetilde{\nu},\widetilde{\lambda})$ with
$p_{\mu,\widetilde{\nu},\widetilde{\lambda}}(M_0)\neq0.$

In this subcase, $w^{'}=\sum
p_{\mu,\widetilde{\nu}}(M_0)M_{-\mu}Y_{-\frac{1}{2}-\widetilde{\nu}}w.$
By using $L_{\underline{k}+1}-\psi(L_{\underline{k}+1})$ to act on
$w^{'},$  we have
\begin{eqnarray*}
&&(L_{\underline{k}+1}-\psi(L_{\underline{k}+1}))w^{'}\\
&=&\sum_{\begin{subarray} \ (\mu,\widetilde{\nu})\in\Lambda_N\\
\ \ \mu(\underline{k})\neq0
\end{subarray}}
p_{\mu,\widetilde{\nu}}(M_0)[L_{\underline{k}+1},
M_{-\mu}]Y_{-\frac{1}{2}-\widetilde{\nu}}w
+\sum_{\begin{subarray}\ (\mu,\widetilde{\nu})\in\Lambda_N\\
\ \ \mu(\underline{k})\neq0
\end{subarray}}
p_{\mu,\widetilde{\nu}}(M_0)
M_{-\mu}[L_{\underline{k}+1},Y_{-\frac{1}{2}-\widetilde{\nu}}]w
\end{eqnarray*}
\begin{equation} \label{eq:1}
+\sum_{\begin{subarray}\ (\mu,\widetilde{\nu})\in\Lambda_N\\
\ \ \mu(\underline{k})=0
\end{subarray}}
p_{\mu,\widetilde{\nu}}(M_0)
[L_{\underline{k}+1},M_{-\mu}Y_{-\frac{1}{2}-\widetilde{\nu}}]w
+\sum_{\begin{subarray}\ (\mu,\widetilde{\nu})\notin\Lambda_N
\end{subarray}}
p_{\mu,\widetilde{\nu}}(M_0)
[L_{\underline{k}+1},M_{-\mu}Y_{-\frac{1}{2}-\widetilde{\nu}}]w.
\end{equation}
We denote the four terms on the right hand side of (3.9) by $v_1,
v_2, v_3$ and $v_4$ respectively. For
$$\mu=(\underline{k}^{\mu(\underline{k})},
(\underline{k}+1)^{\mu(\underline{k}+1)},\cdots),\
\mu(\underline{k})\neq0,$$ we denote
$M_{-\mu}=M_{-\mu^{'}}M_{-\underline{k}}^{\mu(\underline{k})},$
where $\mu^{'}=((\underline{k}+1)^{\mu(\underline{k}+1)},
(\underline{k}+2)^{\mu(\underline{k}+2)},\cdots).$ Note that
$[L_{\underline{k}+1},M_{-\mu^{'}}]\in U(\mathfrak{b}^{-})$ and
$[L_{\underline{k}+1},M_{-\underline{k}}^{\mu(\underline{k})}]
=-\mu(\underline{k})\underline{k}M_{-\underline{k}}^{\mu(\underline{k})-1}M_{1},$
we have
$$v_{1}=v_{1}^{'}-\sum_{\begin{subarray} \
(\mu,\widetilde{\nu})\in\Lambda_N\\
\ \ \mu(\underline{k})\neq0
\end{subarray}}
\mu(\underline{k})\underline{k}\psi(M_1)p_{\mu,\widetilde{\nu}}(M_0)
M_{-\underline{k}}^{\mu(\underline{k})-1}Y_{-\frac{1}{2}-\widetilde{\nu}}w,$$
where $maxdeg(v_{1}^{'})<N-\underline{k}.$ Thus
$maxdeg(v_1)=N-\underline{k}.$ For $v_i$ $(i=2, 3),$ note that
$[L_{\underline{k}+1},Y_{-\frac{1}{2}-\widetilde{\nu}}]\in
U(\mathfrak{b}^{-})$ and
$[L_{\underline{k}+1},M_{-\mu}Y_{-\frac{1}{2}-\widetilde{\nu}}]\in
U(\mathfrak{b}^{-}),$ we have $maxdeg(v_{i})<N-\underline{k}.$
Finally, for $v_4,$ by using Lemma 3.4, we have
$maxdeg(v_{4})<N-\underline{k}.$ Thus
$(L_{\underline{k}+1}-\psi(L_{\underline{k}+1}))w^{'}\neq0.$

\vskip 3mm\noindent{\bf{Subcase 2.}} There exists some
$\widetilde{\lambda}\neq\bar{0}$ for which
$p_{\mu,\widetilde{\nu},\widetilde{\lambda}}(M_0)\neq0.$

Denote
$$N^{'}:=max\{|\mu|+|\frac{1}{2}+\widetilde{\nu}|+|\widetilde{\lambda}|
|\widetilde{\lambda}\neq\bar{0},
p_{\mu,\widetilde{\nu},\widetilde{\lambda}}(M_0)\neq0\},$$ and set
$$\Lambda_{N^{'}}:=\{(\mu,\widetilde{\nu},\widetilde{\lambda})|
\widetilde{\lambda}\neq\bar{0},
p_{\mu,\widetilde{\nu},\widetilde{\lambda}}(M_0)\neq0,
|\mu|+|\frac{1}{2}+\widetilde{\nu}|+|\widetilde{\lambda}|=N^{'}\},$$
\begin{eqnarray*}l:=min\{n|~\widetilde{\lambda}&=&
(n^{\lambda(n)},(n+1)^{\lambda(n+1)},\cdots)\ \mathrm{such}\ \mathrm{that}\\\
&& |\mu|+|\frac{1}{2}+\widetilde{\nu}|+|\widetilde{\lambda}|=N^{'}\
\mathrm{and}\
p_{\mu,\widetilde{\nu},\widetilde{\lambda}}(M_{0})\neq0\}.
\end{eqnarray*}
Note that $\widetilde{\lambda}=\bar{0}$ for those
$(\mu,\widetilde{\nu},\widetilde{\lambda})$ satisfying
$N^{'}<|\mu|+|\frac{1}{2}+\widetilde{\nu}|+|\widetilde{\lambda}|\leq
N$ and $p_{\mu,\widetilde{\nu},\widetilde{\lambda}}(M_0)\neq0.$ Thus
we have
\begin{eqnarray*}
w^{'}&=&\sum_{\begin{subarray}{c}
(\mu,\widetilde{\nu},\widetilde{\lambda})\in\Lambda_{N^{'}}\\
\lambda(l)\neq0
\end{subarray}}
p_{\mu,\widetilde{\nu},\widetilde{\lambda}}(M_0)
M_{-\mu}Y_{-\frac{1}{2}-\widetilde{\nu}}L_{-\widetilde{\lambda}}w\\
&+&\sum_{\begin{subarray}\ (\mu,\widetilde{\nu},\widetilde{\lambda})
\in\Lambda_{N^{'}}\\
\ \ \lambda(l)=0
\end{subarray}}
p_{\mu,\widetilde{\nu},\widetilde{\lambda}}(M_0)
M_{-\mu}Y_{-\frac{1}{2}-\widetilde{\nu}}L_{-\widetilde{\lambda}}w\\
&+&\sum_{\begin{subarray}\
|\mu|+|\frac{1}{2}+\widetilde{\nu}|+|\widetilde{\lambda}|<N^{'}
\end{subarray}}
p_{\mu,\widetilde{\nu},\widetilde{\lambda}}(M_0)
M_{-\mu}Y_{-\frac{1}{2}-\widetilde{\nu}}L_{-\widetilde{\lambda}}w\\
\end{eqnarray*}
\begin{equation}
+\sum_{\begin{subarray}\
N^{'}<|\mu|+|\frac{1}{2}+\widetilde{\nu}|\leq N
\end{subarray}}
p_{\mu,\widetilde{\nu}}(M_0)
M_{-\mu}Y_{-\frac{1}{2}-\widetilde{\nu}}w.
\end{equation}
We apply $(M_{l+1}-\psi(M_{l+1}))$ to act on both sides of (3.10)
and write the resulting four terms on the right hand side by $v_1,$
$v_2,$ $v_3$ and $v_4$ respectively. It is obvious that $v_4=0.$ For
$v_1,$ note that $l>\underline{k}>0,$ by Lemma 3.2 (iii), we see
that
$$v_1=v-\sum_{\begin{subarray} \
(\mu,\widetilde{\nu},\widetilde{\lambda})\in\Lambda_{N^{'}}\\
\ \ \lambda(l)\neq0
\end{subarray}}
p_{\mu,\widetilde{\nu},\widetilde{\lambda}}(M_0)(l+1)\lambda(l)\psi(M_1)
M_{-\mu}Y_{-\frac{1}{2}-\widetilde{\nu}}L_{-\widetilde{\lambda^{'}}}w,$$
where $maxdeg(v)<N^{'}-l$ and $maxdeg(v_1)=N^{'}-l.$ For $v_2,$
since $[M_{l+1},L_{-\widetilde{\lambda}}]\in U(\mathfrak{b}^{-}),$
we have $maxdeg(v_2)\leq N^{'}-l-1.$ Finally for $v_3,$ by using
Lemma 3.2 (i), we see that $maxdeg(v_3)<N^{'}-l.$ These imply that
$(M_{l+1}-\psi(M_{l+1}))w^{'}\neq0.$ The proof of Theorem 3.5 is
completed. \hfill $\Box$

\vskip 3mm\noindent{\bf{Corollary 3.6.}} The center of
$U(\mathfrak{sv})$ is $\mathbb{C}[M_0].$

\vskip 3mm\noindent{\bf{Proof.}} For any $z\in Z(U(\mathfrak{sv})),$
the center of $U(\mathfrak{sv}),$ $zw$ is a Whittaker vector, so
$z\in \mathbb{C}[M_0]$ by Theorem 3.5. This means
$Z(U(\mathfrak{sv}))\subseteq \mathbb{C}[M_0]$ and then
$Z(U(\mathfrak{sv}))=\mathbb{C}[M_0].$

\vskip 3mm\noindent{\bf{Theorem 3.7.}} Suppose $\psi(M_1)\neq0$ and
$\bar{w}=\overline{1\otimes1}\in L_{\psi,\xi}.$ Then $v\in
L_{\psi,\xi}$ is a Whittaker vector if and only if $v=u\bar{w}$ for
some $u\in \mathbb{C}1.$

\vskip 3mm\noindent{\bf{Proof.}} It is easy to see that the set
$$\{M_{-\mu}Y_{-\frac{1}{2}-\widetilde{\nu}}L_{-\widetilde{\lambda}}\bar{w}|
\mu\in \mathcal{P}, \widetilde{\nu}, \widetilde{\lambda}\in
\widetilde{\mathcal{P}}, k\in \mathbb{N}\}$$ forms a basis of
$L_{\psi,\xi}.$ Then we can use the same argument as in Theorem 3.5
to complete the proof of Theorem 3.7. \hfill $\Box$

\vskip 3mm\noindent{\bf{Theorem 3.8.}} Let $\psi_1, \psi_2$ be Lie
algebra homomorphisms from $\mathfrak{sv}^{+}$ to $\mathbb{C}$ and
$\xi_1, \xi_2\in\mathbb{C}.$  Then the $\mathfrak{sv}$-modules
$L_{\psi_1,\xi_1}$ and $L_{\psi_2,\xi_2}$ are isomorphic if and only
if $\psi_1=\psi_2, \xi_1=\xi_2.$

\vskip 3mm\noindent{\bf{Proof.}} Suppose that $\bar{w}_i$ is a
cyclic Whittaker vector of $L_{\psi_i,\xi_i}(i=1,2),$
$f:L_{\psi_1,\xi_1}\rightarrow L_{\psi_2,\xi_2}$ is an isomorphism
of modules. Then
$$E_nf(\bar{w}_1)=f(E_n\bar{w}_1)=\psi_1(E_n)f(\bar{w}_1),\ \forall\
n\in\mathbb{Z}_{+}.$$ Thus $f(\bar{w}_1)$ is a Whittaker vector of
type $\psi_1.$ This implies $\psi_1=\psi_2$ since there are no
Whittaker vectors of type other than $\psi_2$ in $L_{\psi_2,\xi_2}$
by the proof of Theorem 3.5. Moreover,
$\xi_1f(\bar{w}_1)=f(M_0\bar{w}_1)=M_0f(\bar{w}_1)=\xi_2f(\bar{w}_1),$
we get $\xi_1=\xi_2.$ This completes the proof. \hfill $\Box$

\vskip 3mm\noindent{\section{Irreducible Whittaker modules of
nonsingular type}}

\vskip 3mm In this section, the Lie algebra homomorphism $\psi$ is
assumed to be non-singular, that is $\psi(M_1)\neq0.$ we prove that
the Whittaker module $L_{\psi,\xi},$ defined by (2.2), is
irreducible, and we also prove that every irreducible Whittaker
module of type $\psi$ for the Schr\"{o}dinger Virasoro algebra
$\mathfrak{sv}$ is isomorphic to $L_{\psi,\xi}$ for some $\xi\in
\mathbb{C}$.

\vskip 3mm Fix a Whittaker module $V$ of type $\psi$ with cyclic
Whittaker vector $w.$ $V$ is naturally a $\mathfrak{sv}^{+}$-module.
Following [10] and [18] we define a new action, called dot action,
of $\mathfrak{sv}^{+}$ on $V$ by setting
\begin{eqnarray}
x\cdot v=xv-\psi(x)v,\ \mathrm{for}\ x\in \mathfrak{sv}^{+}\
\mathrm{and}\ v\in V.
\end{eqnarray}
Then it is clear that $V$ is a $\mathfrak{sv}^{+}$-module under the
dot action, and we have $E_n\cdot v=E_nv-\psi(E_n)v=[E_n,u]w$ for
$n\in \mathbb{Z}_{+}$ and $v=uw\in V.$

\vskip 3mm\noindent{\bf{Lemma 4.1.}} If $n\in \mathbb{Z}_{+},$ then
$E_n$ acts locally nilpotent on $V$ under the dot action.

\vskip 3mm\noindent{\bf{Proof.}} By the Lie products of
$\mathfrak{sv},$ we see that
$$ad^{2}M_n(L_i)=0=adM_n(M_i)=adM_n(Y_{\frac{1}{2}+i}),
 \forall i\in\mathbb{Z}.$$
$$ad^{3}Y_{\frac{1}{2}+(n-1)}(L_i)=0=adY_{\frac{1}{2}+(n-1)}(M_i)
=ad^{2}Y_{\frac{1}{2}+(n-1)}(Y_{\frac{1}{2}+i}), \forall
i\in\mathbb{Z}.$$ Thus, for any basis element
$u=M_{0}^{k}M_{-\mu}Y_{-\frac{1}{2}-\widetilde{\nu}}L_{-\widetilde{\lambda}}$
of $U(\mathfrak{b}^{-}),$ it is clear that $ad^{2}M_n,$
$ad^{3}Y_{\frac{1}{2}+(n-1)}$ act on $u$ as zero.  To prove
$ad^{m}L_n(u)=0$ for $m$ sufficiently large, we note that
$ad^{m}L_n(M_{-\mu}Y_{-\frac{1}{2}-\widetilde{\nu}}L_{-\widetilde{\lambda}})
\in
U(\mathfrak{sv})_{-(|\mu|+|\frac{1}{2}+\widetilde{\nu}|+|\widetilde{\lambda}|)+nm},$
so
$ad^{m}L_n(M_{-\mu}Y_{-\frac{1}{2}-\widetilde{\nu}}L_{-\widetilde{\lambda}})$
is a combination of basis elements of $U(\mathfrak{sv})$ of the form
\begin{eqnarray}
M_{0}^{k}M_{-\mu_1}Y_{-\frac{1}{2}-\widetilde{\nu}_1}
L_{-\widetilde{\lambda}_1}M_{m_p}\cdots
M_{m_{1}}Y_{\frac{1}{2}+n_q}\cdots
Y_{\frac{1}{2}+n_{1}}L_{l_h}\cdots L_{l_{1}},
\end{eqnarray}
where
$-(|\mu_1|+|\frac{1}{2}+\widetilde{\nu}_{1}|+|\widetilde{\lambda}_1|)
+\sum_{i=1}^{p}m_{i}+\sum_{i=1}^{q}(\frac{1}{2}+n_{i})+\sum_{i=1}^{h}l_{i}
=-(|\widetilde{\lambda}|+|\mu|+|\frac{1}{2}+\widetilde{\nu})|)+nm,$
$\#(\mu_1)+\#(\widetilde{\nu}_1)+\#(\widetilde{\lambda}_1)+p+q+h\leq
\#(\mu,\widetilde{\nu},\widetilde{\lambda}).$ Recall that $E_nw=0$
if $n\geq3,$ it is easy to see that the element, given in (4.2),
acts on the cyclic Whittaker vector $w$ as zero when $m$
sufficiently large. This finishes the proof of the lemma. \hfill
$\Box$

\vskip 3mm\noindent{\bf{Lemma 4.2.}} Let
$(\mu,\widetilde{\nu},\widetilde{\lambda})\in
\mathcal{P}\times\widetilde{\mathcal{P}}\times\widetilde{\mathcal{P}},
k\in \mathbb{N}.$

(i) For all $n>0, E_n\cdot
M_0^{k}M_{-\mu}Y_{-\frac{1}{2}-\widetilde{\nu}}L_{-\widetilde{\lambda}}w
\in \mathrm{span}_{\mathbb{C}}\{M_0^{i}
M_{-\mu^{'}}Y_{-\frac{1}{2}-\widetilde{\nu^{'}}}L_{-\widetilde{\lambda^{'}}}w|\
\mathrm{where}\
|\mu^{'}|+|\frac{1}{2}+\widetilde{\nu^{'}}|+|\widetilde{\lambda^{'}}|
+\lambda^{'}(0)\leq
|\mu|+|\frac{1}{2}+\widetilde{\nu}|+|\widetilde{\lambda}|+\lambda(0);i=k,k+1\}.$

(ii) If
$n>|\mu|+|\frac{1}{2}+\widetilde{\nu}|+|\widetilde{\lambda}|+2,$
then $E_n\cdot
(M_{-\mu}Y_{-\frac{1}{2}-\widetilde{\nu}}L_{-\widetilde{\lambda}}w)=0.$

\vskip 3mm\noindent{\bf{Proof.}} (i) Since $$E_n\cdot
(M_0^{k}M_{-\mu}Y_{-\frac{1}{2}-\widetilde{\nu}}L_{-\widetilde{\lambda}}w)
=M_0^{k}(E_n\cdot
(M_{-\mu}Y_{-\frac{1}{2}-\widetilde{\nu}}L_{-\widetilde{\lambda}}w)),$$
we only need to prove (i) for $k=0.$ The result for
$\#(\mu,\widetilde{\nu},\widetilde{\lambda})=0$ is obvious. Now we
prove the result for $\#(\mu,\widetilde{\nu},\widetilde{\lambda})>0$
by induction.

For the case $\mu\neq \bar{0},$ set $m=max\{i|\mu(i)>0\}.$ Then
$$M_{-\mu}Y_{-\frac{1}{2}-\widetilde{\nu}}L_{-\widetilde{\lambda}}
=M_{-m}M_{-\mu^{'}}Y_{-\frac{1}{2}-\widetilde{\nu}}L_{-\widetilde{\lambda}},$$
where $\mu^{'}(m)=\mu(m)-1, \mu^{'}(i)=\mu(i)$ for all $i\neq k.$
Therefore
\begin{eqnarray*}
E_n\cdot
M_{-\mu}Y_{-\frac{1}{2}-\widetilde{\nu}}L_{-\widetilde{\lambda}}w
\end{eqnarray*}
\begin{eqnarray}
=[E_n,M_{-m}]M_{-\mu^{'}}Y_{-\frac{1}{2}-\widetilde{\nu}}L_{-\widetilde{\lambda}}w
+M_{-m}[E_n,M_{-\mu^{'}}Y_{-\frac{1}{2}-\widetilde{\nu}}L_{-\widetilde{\lambda}}]w.
\end{eqnarray}
For the first term on the right hand side of (4.3), note that
$[E_n,M_{-m}]=0$ for $E_{n}=M_{n}$ or $Y_{\frac{1}{2}+(n-1)},$ we
only need to consider the case for $E_n=L_n.$ If $n-m\leq0,$ it is
obvious
that$$[L_n,M_{-m}]M_{-\mu^{'}}Y_{-\frac{1}{2}-\widetilde{\nu}}
L_{-\widetilde{\lambda}}w=-mM_{n-m}M_{-\mu^{'}}Y_{-\frac{1}{2}-\widetilde{\nu}}
L_{-\widetilde{\lambda}}w$$ has the desired form. If $n-m>0,$
\begin{eqnarray*}
M_{n-m}M_{-\mu^{'}}Y_{-\frac{1}{2}-\widetilde{\nu}}L_{-\widetilde{\lambda}}w
=M_{n-m}\cdot (M_{-\mu^{'}}Y_{-\frac{1}{2}-\widetilde{\nu}}
L_{-\widetilde{\lambda}}w)
+\psi(M_{n-m})M_{-\mu^{'}}Y_{-\frac{1}{2}-\widetilde{\nu}}L_{-\widetilde{\lambda}}w.
\end{eqnarray*}
By assumption,  $M_{n-m}\cdot
(M_{-\mu^{'}}Y_{-\frac{1}{2}-\widetilde{\nu}}
L_{-\widetilde{\lambda}}w),$ and therefore
$M_{n-m}M_{-\mu^{'}}Y_{-\frac{1}{2}-\widetilde{\nu}}L_{-\widetilde{\lambda}}w,$
has the desired form. For the second term on the right hand side of
(4.3), we have, by the induction hypothesis, that
\begin{eqnarray*}
[E_n,M_{-\mu^{'}}Y_{-\frac{1}{2}-\widetilde{\nu}}L_{-\widetilde{\lambda}}]w&\in&
\mathrm{span}_{\mathbb{C}}\{M_{0}^{i}
M_{-\mu^{''}}Y_{-\frac{1}{2}-\widetilde{\nu^{'}}}L_{-\widetilde{\lambda^{'}}}w|
|\mu^{''}|+|\frac{1}{2}+\widetilde{\nu^{'}}|+|\widetilde{\lambda^{'}}|\\
&&+\lambda^{'}(0)\leq
|\mu^{'}|+|\frac{1}{2}+\widetilde{\nu}|+|\widetilde{\lambda}|+\lambda(0);i=0,1;j=0,1\}.
\end{eqnarray*}
Thus
$$M_{-m}[E_n,M_{-\mu^{'}}Y_{-\frac{1}{2}-\widetilde{\nu}}L_{-\widetilde{\lambda}}]w$$
has the desired form since $-m<0$ and
$m+|\widetilde{\lambda^{'}}|+|\mu|+|\frac{1}{2}+\widetilde{\nu}|
=|\widetilde{\lambda}|+|\mu|+|\frac{1}{2}+\widetilde{\nu}|.$

For the case  $\mu=\bar{0},\widetilde{\nu}\neq\bar{0}$ or
$\mu=\widetilde{\nu}=\bar{0}, \widetilde{\lambda}\neq\bar{0},$ one
can prove the result by a similar argument as we did in the first
case. This is omitted for shortness.

(ii) Note that
$$[E_n,M_{-\mu}Y_{-\frac{1}{2}-\widetilde{\nu}}L_{-\widetilde{\lambda}}]
=\sum_{\mu^{1},\widetilde{\nu^{1}},\widetilde{\lambda^{1}},E_m}
p_{\mu^{1},\widetilde{\nu^{1}},\widetilde{\lambda^{1}},E_m}
M_{-\mu^{1}}Y_{-\frac{1}{2}-\widetilde{\nu^{1}}}L_{-\widetilde{\lambda^{1}}}E_m,$$
where
$p_{\mu^{1},\widetilde{\nu^{1}},\widetilde{\lambda^{1}},E_m}\in
\mathbb{C},$ and
$m=n-(|\mu|+|\frac{1}{2}+\widetilde{\nu}|+|\widetilde{\lambda}|)
+(|\mu^{1}|+|\frac{1}{2}+\widetilde{\nu^{1}}|+|\widetilde{\lambda^{1}}|)>2,
$ This implies (ii) as $E_mw=0$ for any $m>2.$\hfill $\Box$

\vskip 5mm\noindent{\bf{Lemma 4.3.}} Suppose $V$ is a Whittaker
module for $\mathfrak{sv},$ and let $v\in V.$ Regarding $V$ as an
$\mathfrak{sv}^{+}$-module under the dot action, then
$U(\mathfrak{sv}^{+})\cdot v$ is a finite-dimensional
$\mathfrak{sv}^{+}$- submodule of $V.$

\vskip 3mm\noindent{\bf{Proof.}} This is a direct result of Lemma
4.2.\hfill $\Box$

\vskip 5mm\noindent{\bf{Lemma 4.4.}} Let $V$ be a Whittaker module
for $\mathfrak{sv},$ and let $S\subseteq V$ be a nonzero submodule.
Then there is a nonzero Whittaker vector $w^{'}\in S.$

\vskip 3mm\noindent{\bf{Proof.}} $\forall 0\neq v\in S,$ by Lemma
4.3, $U(\mathfrak{sv}^{+})\cdot v$ is a finite-dimensional submodule
of $S.$ Then by Lemma 4.1, we know that every element $E_n$ of
$\mathfrak{sv}^{+}$ is nilpotent on $U(\mathfrak{sv}^{+})\cdot v$
under the dot action. Then the result of the lemma follows from
Engel's Theorem (Theorem 3.3 in [8]).\hfill $\Box$

\vskip 5mm\noindent{\bf{Proposition 4.5.}} For any $\xi\in
\mathbb{C},$ the Whittaker module $L_{\psi,\xi}$ for the
Schr\"{o}dinger Virasoro algebra $\mathfrak{sv}$ is irreducible.

\vskip 3mm\noindent{\bf{Proof.}} It follows from Lemma 4.4 and
Theorem 3.7.\hfill $\Box$

\vskip 3mm It is known (see Lemma 2.1.3 in [4]) that the Schur's
Lemma can be generalized to infinite dimensional irreducible modules
with countable cardinality.

\vskip 5mm\noindent{\bf{Theorem 4.6.}} Let $S$ be an irreducible
Whittaker module of type $\psi$ for the Schr\"{o}dinger Virasoro
algebra $\mathfrak{sv}.$ Then $S\cong L_{\psi,\xi}$ for some $\xi\in
\mathbb{C}.$

\vskip 3mm\noindent{\bf{Proof.}} Let $w_s \in S$ be a cyclic
Whittaker vector corresponding to $\psi.$ Since $M_0$ acts by a
scalar by Schur's Lemma, there exists $\xi\in \mathbb{C}$ such that
$M_0s=\xi s$ for all $s\in S.$ Now by the universal property of
$W_{\psi},$ there exists a module homomorphism
$\varphi:W_{\psi}\rightarrow S$ with $uw\mapsto uw_s.$ This map is
surjective since $w_s$ generates $S.$ But then
$$\varphi((M_0-\xi)W_{\psi})=(M_0-\xi)\varphi(W_{\psi})=(M_0-\xi)S=0,$$
so it follows that $$(M_0-\xi)W_{\psi}\subseteq ker\varphi\subseteq
W_{\psi}.$$ Because $L_{\psi,\xi}$ is irreducible by Proposition 4.5
and $ker\varphi\neq W_{\psi},$ this forces
$ker\varphi=(M_0-\xi)W_{\psi}.$\hfill $\Box$

\vskip 3mm For a given $\psi:\mathfrak{sv}^{+}\rightarrow\mathbb{C}$
and $\xi\in\mathbb{C},$ note that
\begin{eqnarray*}
I&=&U(\mathfrak{sv})(M_{0}-\xi)+\sum_{i\in\mathbb{Z}_{+}}U(\mathfrak{sv})
(L_{i}-\psi(L_{i}))+\sum_{i\in\mathbb{Z}_{+}}U(\mathfrak{sv})
(M_{i}-\psi(M_{i}))\\
&+&\sum_{i\in\mathbb{N}}U(\mathfrak{sv})
(Y_{\frac{1}{2}+i}-\psi(Y_{\frac{1}{2}+i}))
\end{eqnarray*}
is a left ideal of $U(\mathfrak{sv}).$ For $u\in U(\mathfrak{sv}),$
let $\bar{u}$ denote the coset $u+I\in U(\mathfrak{sv})/I.$ Then we
may regard $U(\mathfrak{sv})/I$ as a Whittaker module of type $\psi$
with cyclic Whittaker vector $\bar{1}.$ We have the following
result.

\vskip 5mm\noindent{\bf{Lemma 4.7.}} The Whittaker module
$V=U(\mathfrak{sv})/I$ is irreducible, and thus $V\cong
L_{\psi,\xi}.$

\vskip 5mm\noindent{\bf{Theorem 4.8.}} Suppose that $V$ is a
Whittaker module of type $\psi$ such that $M_0$ acts by a scalar
$\xi\in \mathbb{C},$ then $V$ is irreducible.

\vskip 3mm\noindent{\bf{Proof.}} Let $K$ denote the kernel of the
natural surjective map $U(\mathfrak{sv})\rightarrow V$ given by
$u\mapsto uw,$ where $w$ is a cyclic Whittaker vector of $V.$  Then
$K$ is a proper left ideal containing $I.$ By Lemma 4.7, $I$ is
maximal, and thus $K=I$ and $V\cong U(\mathfrak{sv})/I$ is
irreducible.\hfill $\Box$

\vskip 3mm\noindent{\section{Whittaker modules of singular type}}

\vskip 3mm From now on, we assume that the Lie algebra homomorphism
$\psi$ is singular, that is $\psi(M_1)=0.$

\vskip 3mm As in Theorem 3.5 and Theorem 3.7, we use the notation
$w$ (resp. $\bar{w}$) to denote the cyclic Whittaker vector
$1\otimes1$ (resp. $\overline{1\otimes1}$) for $W_{\psi}$ (resp.
$L_{\psi,\xi}$). The following facts about $W_{\psi}$ and
$L_{\psi,\xi}$ are obvious:

\vskip 3mm (i) $W_{\psi}=U(\mathfrak{sv})w$ is free as a
$U(\mathfrak{b}^{-})$-module and the set
\begin{equation} \label{eq:1}
\{M_0^kM_{-\mu}Y_{-\frac{1}{2}-\widetilde{\nu}}L_{-\widetilde{\lambda}}\bar{w}|~
(\mu,\widetilde{\nu},\widetilde{\lambda})\in \mathcal{P}\times
\widetilde{\mathcal{P}}\times\widetilde{\mathcal{P}},k\in
\mathbb{N}\}
\end{equation}
forms a basis for $W_{\psi}$ by the PBW Theorem.

 \vskip 3mm (ii) $L_{\psi,\xi}=U(\mathfrak{sv})\bar{w}$ is
free as a $U(\mathfrak{sv}^{-}\oplus\mathbb{C}L_0)$-module and the
set
\begin{equation} \label{eq:1}
\{M_{-\mu}Y_{-\frac{1}{2}-\widetilde{\nu}}L_{-\widetilde{\lambda}}\bar{w}|
(\mu,\widetilde{\nu},\widetilde{\lambda})\in \mathcal{P}\times
\widetilde{\mathcal{P}}\times\widetilde{\mathcal{P}}\}
\end{equation}
forms a basis for $L_{\psi,\xi}$.

\vskip 3mm (iii) $W_{\psi}$ and $L_{\psi,\xi}$ are
$\mathfrak{sv}^{+}$-modules under the dot action defined in (4.1).

\vskip 3mm Denote $\psi(L_1)=\eta_1, \psi(L_2)=\eta_2$ and
$\psi(Y_\frac{1}{2})=\eta_3.$ Set $$ z=\left\{\begin{array}{cl} L_0
& \textrm{if }\eta_1=\eta_2=\eta_3=0,\\
\!\!\!\!L_0M_0^2\!-\!\eta_2M_{-2}M_0\!-\!\eta_1M_{-1}M_0\!+\!\eta_2M_{-1}^2\!-\!
\frac{\eta_3}{2}Y_{-\frac{1}{2}}M_0\!+\!\frac{\eta_3^2}{2}M_{-1}
&\textrm{if } \eta_2\neq 0 \textrm{ or } \eta_3\neq 0,\\
L_0M_0-\eta_1M_{-1} &\textrm{if   others.}
\end{array}\right.$$

\vskip 5mm\noindent{\bf{Proposition 5.1.}}  For $u\in
\mathbb{C}[z,M_0]$, $v=uw$ is a Whittaker vector of $W_{\psi}$.

\vskip 3mm\noindent{\bf{Proof.}} It is easy to check that $L_1\cdot
v=L_2\cdot v=M_1\cdot v=Y_{\frac{1}{2}}\cdot v=0.$ Thus for any
$y\in \mathfrak{sv}^{+},$ we have $y\cdot v=0.$\hfill $\Box$

\vskip 5mm\noindent{\bf{Proposition 5.2.}}  If we set
 $$ z=\left\{\begin{array}{cl}
L_0 & \textrm{if
}\eta_1=\eta_2=\eta_3=0,\\
\!\!\!\!\xi^2L_0\!-\!\xi\eta_2M_{-2}\!-\!\xi\eta_1M_{-1}\!+\!\eta_2M_{-1}^2\!-\!
\xi\frac{\eta_3}{2}Y_{-\frac{1}{2}}\!+\!\frac{\eta_3^2}{2}M_{-1}
& \textrm{if } \eta_2\neq 0 \textrm{ or } \eta_3\neq 0,\\
\xi L_0-\eta_1M_{-1} &\textrm{if   others,}
\end{array}\right.$$
then for $u\in \mathbb{C}[z]$, $u\bar{w}$ is a Whittaker vector of
$L_{\psi,\xi}$.

\vskip 3mm\noindent{\bf{Proof.}} It is easy to check, we omit the
details.\hfill $\Box$

\vskip 5mm\noindent{\bf{Theorem 5.3.}} If $\psi$ is singular, then
$L_{\psi,\xi}$ is reducible for any $\xi\in\mathbb{C}.$

\vskip 3mm\noindent{\bf{Proof.}} By Proposition 5.2, we can easily
see that the submodule $V$ generated by $z\bar{w}\in L_{\psi,\xi}$
is a proper Whittaker submodule.\hfill $\Box$

\vskip 5mm If $\psi$ is identically zero, that is
$\eta_1=\eta_2=\eta_3=\psi(M_1)=0,$ then we have the following more
precise results.

\vskip 3mm\noindent{\bf{Theorem 5.4.}} If $\psi$ is identically
zero, then the set of Whittaker Vectors of $W_{\psi}$ is
$\mathbb{C}[L_0,M_0]w.$

\vskip 3mm\noindent{\bf{Proof.}} By Proposition 5.1, we see that
each element of $\mathbb{C}[M_0,L_0]w$ is a Whittaker vector of
$W_{\psi}$. For any $v\in W_{\psi}\setminus\mathbb{C}[M_0,L_0]w,$
noting that $W_{\psi}$ has a basis given by (5.1), we can write it
as
\begin{equation} \label{eq:1}
v=\sum a^{\mu,\widetilde{\nu},\lambda}_{k,l}M_{-\mu}
Y_{-\frac{1}{2}-\widetilde{\nu}}L_{-\lambda}L_0^{k}M_{0}^{l}w,\end{equation}
where $\mu, \lambda\in \mathcal{P}$, $\widetilde{\nu}\in
\widetilde{\mathcal{P}}$ satisfying $\mu\neq\bar{0}$ or
$\widetilde{\nu}\neq\bar{0}$ or $\lambda\neq\bar{0}$ for some
$a^{\mu,\widetilde{\nu},\lambda}_{k,l}\neq0.$ For $v$ defined in
(5.3), we need to prove that there exists $x\in \mathfrak{sv}^{+}$
such that $x\cdot v\neq0.$

\vskip 3mm\noindent{\bf{Case 1.}} There exists $\lambda\neq\bar{0}$
such that $a^{\mu,\widetilde{\nu},\lambda}_{k,l}\neq0$ in (5.3).

\vskip 3mm\noindent We denote
$$p=\max\{\#(\lambda)~|~a^{\mu,\widetilde{\nu},\lambda}_{k,l}\neq0\}$$
and
$$N=\max\{\lambda_{p}~|~a^{\mu,\widetilde{\nu},\lambda}_{k,l}\neq0\}.$$
Then $N\geq1$ and $v$ is of the form
\begin{eqnarray}
v=\sum_{\lambda_p=N} a^{\mu,\widetilde{\nu},\lambda}_{k,l}M_{-\mu}
Y_{-\frac{1}{2}-\widetilde{\nu}}L_{-\lambda}L_0^{k}M_0^{l}w +\sum_{
\lambda_p<N } a^{\mu,\widetilde{\nu},\lambda}_{k,l}M_{-\mu}
Y_{-\frac{1}{2}-\widetilde{\nu}}L_{-\lambda}L_0^{k}M_0^{l}w.
\end{eqnarray}
 By using $M_{N}$ to act on both sides of (5.4) by dot
action, we see that
\begin{eqnarray*}M_{N}\cdot
v&=&-\sum_{\#(\lambda^{'})=p-1}\lambda(N)N
a^{\mu,\widetilde{\nu},\lambda}_{k,l}M_{-\mu}
Y_{-\frac{1}{2}-\widetilde{\nu}}L_{-\lambda^{'}}L_0^{k}M_0^{l+1}w
\\&&+\sum_{\#(\lambda^{''})<p-1}
b^{\mu,\widetilde{\nu},\lambda^{''}}_{k,l}M_{-\mu}
Y_{-\frac{1}{2}-\widetilde{\nu}}L_{-\lambda^{''}}L_0^{k}M_0^{l}w\\&\neq&0,
\end{eqnarray*} where $\lambda^{'}(i)=\lambda(i)$ if $i\neq p$ and
$\lambda^{'}(p)=\lambda(p)-1,$
$b^{\mu,\widetilde{\nu},\lambda^{''}}_{k,l}\in\mathbb{C}.$

\vskip 3mm\noindent{\bf{Case 2.}} $\lambda=\bar{0}$ for any
$a^{\mu,\widetilde{\nu},\lambda}_{k,l}\neq0$  and there exists
$\widetilde{\nu}\neq\bar{0}$ such that
$a^{\mu,\widetilde{\nu},\lambda}_{k,l}\neq0$ in (5.3).

\vskip 3mm\noindent In this case, (5.3) becomes
\begin{equation} \label{eq:1}
v=\sum a^{\mu,\widetilde{\nu}}_{k,l}M_{-\mu}
Y_{-\frac{1}{2}-\widetilde{\nu}}L_{0}^{k}M_{0}^{l}w.
\end{equation}
For (5.5), we set
$$b:=max\{\nu_s~|~\widetilde{\nu}=(\nu_1,\cdots,\nu_s),
a^{\mu,\widetilde{\nu}}_{k,l}\neq0\}.$$ Then (5.5) can be rewritten
as
\begin{equation}
v=\sum_{ \nu(b)\neq0 } a^{\mu,\widetilde{\nu}}_{k,l}M_{-\mu}
Y_{-\frac{1}{2}-\widetilde{\nu}}L_{0}^{k}M_{0}^{l}w +\sum_{\nu(b)=0}
a^{\mu,\widetilde{\nu}}_{k,l}M_{-\mu}
Y_{-\frac{1}{2}-\widetilde{\nu}}L_{0}^{k}M_{0}^{l}w.\end{equation}
By using $Y_{\frac{1}{2}+b}$ to act on both sides of (5.6) by dot
action, we have
\begin{eqnarray*}Y_{\frac{1}{2}+b}\cdot v=\sum_{\nu(b)\neq0}
(-1-2b)\nu(b)a^{\mu,\widetilde{\nu}}_{k,l}M_{-\mu}
Y_{-\frac{1}{2}-\widetilde{\nu}^{'}}
Y_{-\frac{1}{2}-b}^{\nu(b)-1}L_{0}^{k}M_{0}^{l+1}w\neq 0,
\end{eqnarray*}
where $\widetilde{\nu}^{'}(i)=\widetilde{\nu}(i)$ for $i\neq b$ and
$\widetilde{\nu}^{'}(b)=0.$

\vskip 3mm\noindent{\bf{Case 3.}} $\lambda=\bar{0}=\widetilde{\nu}$
for any $a^{\mu,\widetilde{\nu},\lambda}_{k,l}\neq0$  in (5.3).

\vskip 3mm\noindent In this case, (5.3) becomes
\begin{equation} \label{eq:1}
v=\sum a^{\mu}_{k,l}M_{-\mu}L_{0}^{k}M_{0}^{l}w.
\end{equation}
For (5.7), we set
$$c:=max\{\mu_{t}~|~\mu=(\mu_1,\cdots,\mu_t), a^{\mu}_{k,l}\neq0\}.$$
Then (5.7) can be rewritten as
\begin{equation} \label{eq:1}
v=\sum_{\mu(c)\neq0}
a^{\mu}_{k,l}M_{-\mu}L_{0}^{k}M_{0}^{l}w+\sum_{\mu(c)=0}a^{\mu}_{k,l}M_{-\mu}L_{0}^{k}M_{0}^{l}w.
\end{equation}
By using $L_{c}$ to act on both sides of (5.8) by dot action,  we
have
$$L_{c}\cdot v=\sum_{\mu(c)\neq0} (-c)\mu(c)a^{\mu}_{k,l}M_{-\mu^{'}}M_{-c}^{\mu(c)-1}L_{0}^{k}M_{0}^{l+1}w
\neq0,$$ where $\mu^{'}(i)=\mu^{i}$ for $i\neq a$ and
$\mu^{'}(a)=0.$ This completes the proof.\hfill $\Box$

\vskip 3mm\noindent{\bf{Theorem 5.5.}} If $\psi$ is identically
zero, $\xi\neq0,$ then the set of Whittaker vectors of
$L_{\psi,\xi}$ is $\mathbb{C}[L_{0}]\bar{w}.$

\vskip 3mm\noindent{\bf{Proof.}} Noting that $L_{\psi,\xi}$ has a
basis defined by (5.2), we can repeat the proof of Theorem 5.4 word
for word except that $M_0$ should be replaced by $\xi$ and $w$
replaced by $\bar{w}.$\hfill $\Box$

\vskip 3mm  Recall the definition of Verma module of generalized
Schr\"{o}dinger-Virasoro algebras given in [21]. We observe that if
$\psi$ is identically zero and $V_{\zeta}$ is the submodule of
$L_{\psi,\xi}$ generated by $(L_0-\zeta)\bar{w},$ where
$\zeta\in\mathbb{C},$ then the quotient module
$$V(\xi,\zeta):=L_{\psi,\xi}/V_{\zeta}$$
is the Verma module for $\mathfrak{sv}.$ Denote by $\bar{\bar{w}}$
the homomorphic image of $\bar{w}$, we immediately obtain the
following Lemma by Theorem 4.6 of [21]:

\vskip 3mm\noindent{\bf{Lemma 5.6.}} The Verma module $V(\xi,\zeta)$
is irreducible if and only if $\xi\neq 0.$

\vskip 3mm\noindent{\bf{Theorem 5.7.}} If $\psi$ is identically
zero, $\xi\neq0,$ then

\vskip 3mm (i) For each $\zeta\in \mathbb{C},$ the Whittaker module
$L_{\psi,\xi}$ has the following filtration
$$L_{\psi,\xi}=V^0\supseteq V^1\supseteq\cdots\supseteq V^{i}\supseteq \cdots$$
where $V^{i}$ is a Whittaker submodule of $L_{\psi,\xi}$ defined by
$V^i=U(\mathfrak{sv})(L_{0}-\zeta)^{i}w,$ and $V^{i+1}$ is a maximal
submodule of $V^i.$ More precisely, $V^i/V^{i+1}$ is isomorphic to
the Verma module $V(\xi,\zeta)$.

\vskip 3mm (ii) $L_{\psi,\xi}$ is isomorphic to $V^{i}$ as
$\mathfrak{sv}$-modules for each $i\in \mathbb{N}.$

\vskip 3mm\noindent{\bf{Proof.}} For (i), it is obvious that
$V^{i}/V^{i+1}\backsimeq V(\xi,\zeta)$ according to the definitions
of $V^{i}, i\in\mathbb{N}.$ Then $V^{i+1}$ is a maximal submodule of
$V^{i}$ by Lemma 5.6. Thus (ii) holds.

For (ii), since $\psi$ is identically zero, we can easily check that
the linear map
\begin{eqnarray*}
f  : & L_{\psi,\xi} & \rightarrow V^{i} \\
     & u\bar{w} & \mapsto u(L_0-\zeta)^{i}\bar{w},
\end{eqnarray*}
where $u\in U(\mathfrak{sv}^{-}\oplus\mathbb{C}L_0),$ is an
isomorphism of modules. \hfill $\Box$

\vskip 3mm\noindent{\bf{Proposition 5.8.}} If $\psi$ is identically
zero, $\xi=0,$ then the submodule $V$ of $L_{\psi,0}$ generated by
$L_{-2}\bar{w}$ is a maximal proper submodule. Moreover,
$L_{\psi,0}/V$ is a one-dimensional trivial module.

\vskip 3mm\noindent{\bf{Proof.}} Note that $L_{-i}\bar{w},
M_{-i-1}\bar{w}, Y_{-\frac{1}{2}-i}\bar{w}\in V$ for all
$i\in\mathbb{N}.$ Thus
$M_{-\mu}Y_{-\frac{1}{2}-\widetilde{\nu}}L_{-\widetilde{\lambda}}\bar{w}\in
V$ for all $(\mu,\widetilde{\nu}, \widetilde{\lambda})\in
\mathcal{P}\times\widetilde{\mathcal{P}}\times\widetilde{\mathcal{P}}$
with $\#(\mu,\widetilde{\nu}, \widetilde{\lambda})>0.$ Since
$M_0\bar{w}=\mathfrak{sv}^{+}\bar{w}=0,$ we see that each element of
$V$ is a linear combination of elements with form
$M_{-\mu}Y_{-\frac{1}{2}-\widetilde{\nu}}L_{-\widetilde{\lambda}}
\bar{w}$, $\#(\mu,\widetilde{\nu}, \widetilde{\lambda})>0.$ Thus
$\bar{w}\notin V.$ So $L_{\psi,0}/V$ is a one-dimensional trivial
quotient module and $V$ is a maximal proper submodule of
$L_{\psi,0}$. \hfill $\Box$

\end{document}